\documentclass[preprint]{elsarticle}

\usepackage{lineno}

\usepackage{natbib}
\usepackage{amssymb}
\usepackage{amsthm}

\usepackage{amsmath}

\RequirePackage{fix-cm}
\usepackage[utf8,latin1]{inputenc}
\usepackage[T1]{fontenc}
\usepackage[caption=false]{subfig}

\usepackage{eucal}

\modulolinenumbers[5]

\journal{Journal of Computational Physics}

\begin{document}

\begin{frontmatter}

\title{Time-symmetry, symplecticity and stability of Euler-Maclaurin and Lanczos-Dyche integration}


\author[label4,label0,label1]{Charalampos M. Markakis\corref{mycorrespondingauthor}}
\cortext[mycorrespondingauthor]{Corresponding author}
\ead{c.markakis@damtp.cam.ac.uk}
\author[label2,label1]{Michael F. O'Boyle}
\author[label2,label1]{Derek Glennon}
\author[label1,label3,label5]{Khoa Tran}
\author[label2,label1,label6]{Pablo Brubeck}
\author[label1]{Roland Haas}
\author[label7]{Hsi-Yu Schive}
\author[label8]{K\=oji Ury\=u}

\address[label0]{DAMTP, Centre for Mathematical Sciences, University of Cambridge, Cambridge, UK}
\address[label4]{School of Mathematical Sciences, Queen Mary University of London,
London, UK}
\address[label2]{Department of Physics, University of Illinois at Urbana-Champaign,
Urbana, IL,
USA}
\address[label1]{NCSA, University of Illinois at Urbana-Champaign, Urbana, IL,
USA}
\address[label3]{Department of Mathematics, University of Illinois at Urbana-Champaign, Urbana, IL, USA}
\address[label5]{Department of Mathematics, University of California San Diego, San Diego, CA, USA}
\address[label6]{Argonne National Laboratory, Lemont, IL, USA}
\address[label7]{Physics Department \& Institute of Astrophysics, National Taiwan University, Taiwan}
\address[label8]{Department of Physics, University of the Ryukyus, Senbaru, Nishihara, Okinawa, Japan}


\begin{abstract}
Numerical evolution of time-dependent differential equations  via explicit Runge--Kutta or Taylor methods  typically fails to preserve symmetries of a system. It is known that there  exists no numerical integration method that in general preserves both the energy and the symplectic structure of a Hamiltonian system. One is thus normally forced to make a choice. Nevertheless, a symmetric integration formula, obtained by Lanczos-Dyche via two-point
Taylor expansion (or Hermite interpolation), is shown here to preserve
both energy as well as symplectic structure for linear
systems. This  formula shares similarities  with the Euler-Maclaurin formula,
but is superconvergent rather than asymptotically convergent.   For partial differential equations,   the  resulting evolution methods are unconditionally stable, i.e, not subject to a 
  Courant-Friedrichs-Lewy limit.  
Although generally implicit, these methods 
become explicit for linear systems. 
\end{abstract}


\end{frontmatter}


\section{Introduction}

In this paper we consider  the numerical solution of initial value problems of the form 
\begin{equation} \label{eq:dudt}
\frac{du}{dt}=f(t,u), \qquad u(t_0)=u_0
\end{equation}
with  $f\!:\Omega \subset \mathbb{R}\times  \mathbb{R}^{N} \rightarrow\mathbb{R}^{N}
$ where $\Omega$ is an open set of $ \mathbb{R}\times  \mathbb{R}^{N} $.
Eq.~(\ref{eq:dudt}) may represent an ordinary differential equation (ODE) for the scalar function $u(t)$, or a system of ODEs, or (in a method of lines approach 
\cite{alcubierreIntroductionNumericalRelativity2012}) a  partial differential equation (PDE) for the vector $\mathbf u(t)$ with components  $u_i(t)=u(t,x_i),$
 constructed from discretizing a  scalar field $u(t,x)$ on the
spatial grid-points $\{x_i\}$.

To evolve this differential equation numerically, one may use the fundamental theorem of calculus to write it as an
integral equation
\begin{equation} \label{eq:ustepintf}
u(t_{\nu+1})=u(t_{\nu})+\int_{t_\nu}^{t_{\nu+1}}f(t,u(t))dt
\end{equation}
(where $\nu \in \mathbb{N}$ counts the sequential time steps of  numerical integration) and compute the time integral via  a suitable quadrature or other numerical
integration rule \cite{miyatakeCharacterizationEnergyPreservingMethods2016}. Explicit integration schemes,  such as  Runge--Kutta methods, or Taylor
series of the integrand expanded about $t=t_\nu$,  yield numerical integration formulas that violate symmetry under time reversal  $t_\nu \leftrightarrow t_{\nu+1} $. Symmetries of the continuum system, such as conservation of energy or symplectic structure, are thus generally violated upon discretization.

It
is known that there  exists no numerical integration method that in
general preserves both the energy and the  symplectic structure of a Hamiltonian
system \cite{butcherSymmetricGeneralLinear2016,chartierAlgebraicApproachInvariant2006,zhongLiePoissonHamiltonJacobiTheory1988}. One is thus normally forced to make a choice between conserving one or the other. However, for quadratic Hamiltonians, we show that it is in fact possible to      preserve
both energy as well as symplectic structure (up to numerical roundoff error). 

Our focus in this paper is to  integrate via a  single-step rule, originally due to Lanczos \cite{lanczosAppliedAnalysis1956}
and Dyche \cite{dycheMultiplePointTaylor1956},
 which is based on  a \textit{two-point} Taylor expansion   about $t=t_\nu$ \textit{and} $t=t_{\nu+1}$,  and is symmetric under time reversal   $t_\nu \leftrightarrow t_{\nu+1} $. We will refer to this integration rule as the Lanczos-Dyche (LD) formula, and will study the degree to which it preserves symmetries of the
continuum system.
  When compared with
Runge--Kutta methods of the
same order,  these symmetric  methods are shown to be efficient, and invariants
of the motion are well--preserved over long
time intervals, exactly for linear systems and approximately for non-linear systems. 

Moreover, when using symmetric methods in a method of lines framework to evolve partial differential equations, the resulting scheme is shown to be  unconditionally stable, that is, not subject to a 
  Courant-Friedrichs-Lewy (CFL) condition, thus allowing arbitrarily  large time
steps in numerical evolution. In  terms of accuracy,  time steps cannot be arbitrarily large, but the high order  of the  methods, combined with the absence of  a CFL limit,  does allow significantly larger time steps compared to explicit methods.  Moreover, while symmetric formulas are implicit, for problems that allow inverting and storing the evolution operator in memory,  there is no computational overhead compared to explicit methods,  as the symmetric formulas can become effectively explicit.

\label{sec1}
\section{Numerical integration via Hermite interpolation}
\label{sec2}

Let $f\!:\mathbb{R}\rightarrow\mathbb{R}$ be a $C^{k}$ function. We wish
to construct  approximations to the definite integral
\begin{equation} \label{eq:DI}
{I(f)}= \int_{{t_1}}^{{t_2}} {f(t)dt}, 
\end{equation}
for a small  interval $\Delta t = t_2-t_1$, given the function values and
derivatives at the endpoints of this interval.

\subsection{One-point Taylor expansion} 

A one-point Taylor expansion of the integrand about $t=t_1$ has the form
\begin{eqnarray} \label{eq:fPRTaylor}
f(t)=\mathcal{P}_n(t)+\mathcal{R}_n(t)
\end{eqnarray}
where
\begin{eqnarray}
\mathcal{P}_n(t)=\sum_{l=0}^{n} \frac{(t-t_1)^l }{l!}f^{(l)}_1
\end{eqnarray}
is the Taylor polynomial of order $n<k$, with $f^{(l)}_1=f^{(l)}(t_1)$,    and
\begin{equation}
\mathcal{R}_n(t) = \frac{(t-t_1)^{n+1}}{(n+1)!} f^{(n+1)}(\tau)
\end{equation}
is the remainder, with
$\tau \in [t_1,t_2]$. Integration of Eq.~\eqref{eq:fPRTaylor} yields the  familiar formula
\begin{equation} \label{eq:Taylorint}
I(f)=
\sum_{l=1}^n \frac{\Delta t^l }{l!} f^{(l-1)}_1 
+R_n(f)
\end{equation}
with remainder
\begin{equation} \label{eq:TaylorintRem}
R_n(f)=\frac{{{\Delta
{t^{n+1}}{}}}}{{(n+1){{!}}}}f^{(n)}(\tau)
\end{equation}
%
%
Eq.~\eqref{eq:dudt} can   be  used to express
all  derivatives $f'_1, f''_1,...$ in terms of $f_1$ 
\cite{butcherNumericalMethodsOrdinary2000}.
This procedure can generally be automated using algebraic manipulation software.
\cite{brownCharacteristicsImplicitMultistep1977,jeltschStabilityPropertiesBrown1978}.
In particular, if Eq.~\eqref{eq:dudt} is a linear ODE,
then  formula~\eqref{eq:Taylorint}
 is  equivalent to a Runge--Kutta method of order $n$.

The above formula has the advantage of being explicit, that is, it only requires initial conditions at time $t=t_1$ to be specified. It is thus suitable for initial value problems such as Eq.~\eqref{eq:dudt}. However, this property means that the  formula is not symmetric under time reversal  $t_1 \leftrightarrow t_2$ and thus fails to preserve  symmetries of the  differential equation \eqref{eq:dudt}.
In what follows, we use two-point Hermite interpolation or, equivalently, a two-point Taylor expansion, to derive formulas symmetric under  time reversal.

\subsection{Two-point Taylor expansion} 
\paragraph{Trapezoidal  rule (2$^{\rm nd}$ order) }
Given the function values $f_1=f(t_1)$ and $f_2=f(t_2)$ at the endpoints,
 one may approximate $f(t)$ in the interval $t\in[t_1,t_2]$ with a     linear
Lagrange interpolating
polynomial,

\begin{equation} 
\mathcal{P}(t)=\frac{{{t_2} - t}}{{\Delta t}}{f_1} + \frac{{t - {t_1}}}{{\Delta
t}}{f_2},
\end{equation}
uniquely determined  by the two conditions  
$\mathcal{P}(t_1)=f_1$ and $\mathcal{P}(t_2)=f_2$. Substituting this expression
into Eq.~\eqref{eq:DI}  and integrating yields the \textit{trapezoidal rule}:
\begin{equation} \label{eq:trapezoidal}
\boxed{ 
{I(f)} =\frac{{\Delta t}}{2}({f_1} + {f_2})+\mathcal{O}(\Delta t^3)
}
\end{equation}
with the error term to be computed below.
\paragraph{Hermite's rule (4$^{\rm th}$ order)}
If, in addition to the function values, the   derivatives $f'^{}_1= f'(t_1)$
and 
 $ f'^{}_2= f'^{}(t_2) $    are  known at the end-points, one may approximate
$f(t)$ with a cubic Hermite interpolating polynomial $\mathcal{P}(t)$,
uniquely determined  by the four conditions  
$\mathcal{P}(t_1)=f_1$,  $\mathcal{P}(t_2)=f_2$,
$\mathcal{P} '(t_1)=f'_1$ and $\mathcal{P} '(t_2)=f'_2$. Using the method
of undetermined coefficients
 \cite{markakisHighorderDifferencePseudospectral2014,nozawaConstructionHighlyAccurate1998}, one finds:
\begin{eqnarray}
\mathcal{P}(t)&=&\frac{{{{({t_2}-t)}^2}(2t + {t_2} - 3{t_1})}}{{\Delta {t^{\rm{3}}}}}{f_1}
-\frac{{(2t + {t_1} - 3{t_2}){{(t - {t_1})}^2}}}{{\Delta {t^{\rm{3}}}}}{f_2}
\nonumber
\\ &+&\frac{{{{({t_2} - t)}^2}(t - {t_1})}}{{\Delta {t^{\rm{2}}}}}{{
f}'_1}
 - \frac{{({t_2} - t){{(t - {t_1})}^2}}}{{\Delta {t^{\rm{2}}}}}{f'^{}_2}.
\end{eqnarray}
 Substituting this expression into Eq.~\eqref{eq:DI}
 and integrating yields the \textit{Hermite rule}:
\begin{equation}  \label{eq:Hermite}
\boxed{
 {I(f)} =\frac{{\Delta t}}{2}({f_1} + {f_2})+\frac{{\Delta t^{2}}}{12}({f'_1}
- {f'_2})+\mathcal{O}(\Delta t^5)
}
\end{equation}
A detailed comparison between Hermite's rule and Simpson's rule may be found in \cite{lampretInvitationHermiteIntegration2004}.

\paragraph{Lotkin's rule (6$^{\rm th}$ order)}

If, in addition to the function values and first derivatives, the  second
derivatives $f''^{}_1=
f''(t_1)$ and 
 $ f''^{}_2= f''^{}(t_2) $   at the end-points are  known, one may approximate
$f(t)$ with a quintic Hermite interpolating polynomial, uniquely determined  by the six conditions  
$\mathcal{P}(t_1)=f_1$,  $\mathcal{P}(t_2)=f_2$,
$\mathcal{P} '(t_1)=f'_1$,  $\mathcal{P}'(t_2)=f'_2$, $\mathcal{P}''(t_1)=f_1$
 and $\mathcal{P}''(t_2)=f_2$.
The method of undetermined coefficients yields:
\begin{eqnarray}
\mathcal{P}(t)&=&\frac{\left[6 t^2+3 t_2 t+10
   t_1^2+t_2^2-5 t_1 \left(3 t+t_2\right)\right]\left(t_2-t\right){}^3 }{\Delta t^{5}}f_1 \nonumber
\\&+&\frac{ \left[6 t^2+3 t_1 t+t_1^2+10
   t_2^2-5 \left(3 t+t_1\right) t_2\right] \left(t-t_1\right){}^3}{\Delta t^{5}}f_2 \nonumber
   \\
   &+&\frac{\left(3 t-4 t_1+t_2\right) \left(t-t_1\right)\left(t_2-t\right){}^3  }{\Delta t^{4}}f'_1 \nonumber
   \\ & +&\frac{\left(3 t+t_1-4 t_2\right)  \left(t-t_1\right){}^3 \left(t_2-t\right)
}{\Delta t^{4}}f'_2 \nonumber
\\&+&\frac{\left(t_2-t\right){}^2 \left(t-t_1\right){}^3 }{2\Delta t^{3}}f''_2+\frac{\left(t_2-t\right){}^3 \left(t-t_1\right){}^2 }{2\Delta t^{3}}f''_1.
\end{eqnarray}
  Substituting this expression into
Eq.~\eqref{eq:DI}
 and integrating yields a $6^{th}$ order  generalization to the trapezoidal rule, which we will refer to as \textit{Lotkin's} \cite{lotkinNewIntegratingProcedure1952} \textit{rule}:\begin{equation} \label{eq:Lotkin}
\boxed{
I(f) = \frac{\Delta t}{2}({f_1} + {f_2}) + \frac{{{\Delta t^2}}}{{10}}({f'_1}
- {f'_2}) + \frac{{{\Delta t^3}}}{{120}}({f''_1} + {f''_2})+\mathcal{O}(\Delta
t^7)
}
\end{equation}
 Unlike formulas such  as Eq.~\eqref{eq:Taylorint} stemming from one-point
Taylor series, one  notices that  Eqs.~\eqref{eq:Hermite} and \eqref{eq:Lotkin} have  different coefficients, which   depend on the order of the series truncation. This is  a characteristic of two-point Taylor series, as shown below.

\paragraph{The Lanczos-Dyche  formula ($2n^{\rm th}$ order)}

The above procedure can be continued up to $2n^{\rm th}$ order if the function values and  $l^{\rm th}$
 derivatives $f^{(l)}_1=f^{(l)}(t_1)$, $f^{(l)}_2=f^{(l)}(t_2)$ at the end-points are known for  $l=0,1,...,n-1$. The finite form of the two-point Taylor series is \cite{hummelGeneralizationTaylorExpansion1949,dycheMultiplePointTaylor1956}
\begin{eqnarray} \label{eq:fPR}
f(t)=\mathcal{P}_n(t)+\mathcal{R}_n(t)
\end{eqnarray}
where
\begin{eqnarray}
\mathcal{P}_n(t)=\sum_{l=0}^{n-1}\Delta t^{-(2l+1)} [(t-t_1)(t-t_2)]^l [(t-t_2)a_l +(t-t_1)b_l ]
\end{eqnarray}
is the two-point Hermite interpolation polynomial,
\begin{eqnarray}
{a_l} = {\left[\frac{{{d^l}}}{{d{t^l}}} {\frac{{f(t)}}{{{{(t - {t_2})}^l}}}} \right]_{t = {t_1}}},\quad {b_l} = {\left[ \frac{{{d^l}}}{{d{t^l}}} {\frac{{f(t)}}{{{{(t -t{_1})}^l}}}} \right]_{t = {t_2}}}
\end{eqnarray}
are constant coefficients, and
\begin{eqnarray}
\mathcal{R}_n(t)=\frac{{[(t - {t_1})(t - {t_2})]^n}}{{(2n)!}}{}{f^{(2n)}}(\tau ), \quad \tau  \in [t_1,t_2]
\end{eqnarray}
is the remainder \cite{estesTwopointTaylorSeries1966}.
 Substituting Eq.~\eqref{eq:fPR}  into
Eq.~\eqref{eq:DI}
 and integrating yields the $2n^{\rm th}$ order generalization of the trapezoidal rule, which we will refer to as the \textit{Lanczos-Dyche}\footnote{We derived this formula independently, but  discovered Lanczos' book  and Dyche's thesis when  this work was near completion, and thus named it after these authors. Despite its remarkable accuracy and symmetry properties,  little attention has been paid to  the LD formula in the literature, and we found
no references other than the ones above. }
\cite{lanczosAppliedAnalysis1956,dycheMultiplePointTaylor1956}
(cf. also \cite{squireApplicationsQuadratureDifferentiation1961,squireCommentApproximateCalculation1962,gouldMaclaurinSecondFormula1963,lambertUseHigherDerivatives1963})
\textit{formula}:
\begin{equation} \label{eq:LD}
\boxed{
I(f) = \sum\limits_{l = 1}^n {C_{ln}  \frac{\Delta t^l }{l!}[{f_1^{(l - 1)}} + {{( - 1)}^{l - 1}}{f_2^{(l - 1)}}]}  + R_n(f)
}
\end{equation}
with coefficients 
\begin{equation} \label{eq:Cln}
C_{ln}:=\frac{n!(2n-l)!}{(2n)!(n-l)!},
\end{equation}
and remainder
\begin{equation} \label{eq:LDerror}
R_n(f)= (-1)^n\frac{n!^2}{(2n+1)!(2n)!}\Delta{t^{2n+1}}{f^{(2n)}}(\tau  ).
\end{equation}
Eqs.~\eqref{eq:trapezoidal},~\eqref{eq:Hermite} and \eqref{eq:Lotkin} are special cases of the Lanczos-Dyche (LD) formula \eqref{eq:LD}, which forms the basis of this paper.

Higher-order rules easily follow from the above  formula. For $n=4$, the  LD formula \eqref{eq:LD} gives the $8^{\rm th}$ order rule
\begin{eqnarray} \label{eq:LD4}
I(f) \!  \! &=&  \! \! \frac{\Delta t}{2}({f_1} + {f_2}) + \frac{{{3\Delta t^2}}}{{28}}({f'_1}
- {f'_2}) + \frac{{{\Delta t^3}}}{{84}}({f''_1} + {f''_2}) \nonumber
\\ \! \! &+&  \! \! \frac{{{\Delta t^4}}}{{1680}}({f'''_1}-{f'''_2})+\mathcal{O}(\Delta t^9).
\end{eqnarray}
For $n=5$, the same formula gives the $10^{\rm th}$ order rule
\begin{eqnarray} \label{eq:LD5}
I(f) \! \! &=& \! \! \frac{\Delta t}{2}({f_1} + {f_2}) + \frac{{{\Delta t^2}}}{{9}}({f'_1}
- {f'_2}) + \frac{{{\Delta t^3}}}{{72}}({f''_1} + {f''_2}) \nonumber
\\ \! \! &+&  \! \! \frac{{{\Delta t^4}}}{{1008}}({f'''_1}-{f'''_2}) +\frac{{{\Delta t^5}}}{{30240}}({f^{(4)}_1}+{f^{(4)}_2}) +\mathcal{O}(\Delta
t^{11}).
\end{eqnarray}
and so forth.

As mentioned earlier, the expansion
coefficients of the LD formula \eqref{eq:LD}  depend
 on the order   $n$ of the series truncation.
This
is in contrast  to  familiar  formulas   based on  one-point Taylor series [cf.
Eq.~\eqref{eq:Taylorint}] or the  Euler-Maclaurin formula [cf.~Eq.~\eqref{eq:EM} below]
which have fixed coefficients. This dependence on the order of truncation
is a characteristic of two-point Taylor series, and is partly responsible for their rapid convergence.

In fact, one of the most interesting features of the LD formula \eqref{eq:LD} is its remainder term, given by Eq.~\eqref{eq:LDerror}.  
Inspection of this remainder
term, as well as numerical applications outlined in the next sections, reveal that the LD formula  is  \textit{superconvergent}: when the function derivatives $f^{(l)}(t)$ of order $l=0,1,...,n-1$ are known,
the  formula 
  is accurate to order $2n$ (rather than order $n$, like  the usual Taylor formula \eqref{eq:Taylorint}). 

\subsection{Comparison to the  Euler-Maclaurin formula}   
\label{sec:23}

It is worth comparing the Lanczos-Dyche formula 
\eqref{eq:LD} to the more widely known (albeit often asymptotically convergent) \textit{Euler-Maclaurin (EM) formula}:
\begin{eqnarray} \label{eq:EM}
I(f) = \sum\limits_{l = 1}^n B_{l} {\frac{ \Delta
{t^l}}{l!}[{f_1^{(l - 1)}} + {{( - 1)}^{l - 1}}{f_2^{(l - 1)}}]}  +  R_n(f)
\end{eqnarray}
where, for even $n$, the remainder term is
\cite{eppersonIntroductionNumericalMethods2013}
\begin{equation} \label{eq:EMerror}
R_n(f)=- B_{n+2} \frac{\Delta{t^{n+3}}}{(n+2)!}{f^{(n+2)}}(\tau  ),
\quad \tau  \in (t_1,t_2).
\end{equation}
and $B_l$ denotes the Bernoulli numbers
\begin{equation} \label{eq:Bnumbers}
B_1=\frac{1}{2}, \, B_2=\frac{1}{6}, \, B_3=0, \, B_4=-\frac{1}{30}, \, B_5=0, \, B_6=\frac{1}{42}, ...
\end{equation}
For $n=1$ and $n=2$, the EM formula \eqref{eq:EM} coincides  with the trapezoidal rule  \eqref{eq:trapezoidal}
and Hermite rule \eqref{eq:Hermite} respectively. However, the EM formula differs from the LD formula  \eqref{eq:LD} at higher order. For $n=4$, Eq.~\eqref{eq:EM} yields the $6^{\rm th}$ order rule:
\begin{equation}
I(f) = \frac{\Delta t}{2}(f_1 + f_2) + \frac{\Delta t^2}{12}(f'_1 - f'_2)
- \frac{\Delta t^4}{720}(f'''_1 - f'''_2) + {\cal O}(\Delta{t^{7}})
\end{equation}
which differs from the Lanczos-Dyche  $8^{\rm th}$ order rule
\eqref{eq:LD4}, despite the fact that both formulas use up to third-order derivatives of $f(t)$. 
 
We note that, according to Eq.~\eqref{eq:Bnumbers}, only odd-order derivatives appear in the EM formula \eqref{eq:EM}, and  these derivatives appear with opposite sign at the end-points. As a result, the EM formula can be used to construct a composite (multi-step) rule by summing up the contributions from a sequence of points ($t_1,t_2,...,t_N$). For the composite EM formula, all derivative contributions cancel out, except for the values at the first and last point ($t_1$ and $t_N$). The EM formula is thus often used to convert a sum $\sum_{i}^{} f(t_i)$ to an integral $\int f(t)dt$ \cite{taoCompactnessContradiction2014,taoEulerMaclaurinFormulaBernoulli2010}.
 An important application is  Riemann zeta-function regularization
\cite{hawkingZetaFunctionRegularization1977,elizaldeExpressionsZetaFunction1989,fermiLocalZetaRegularization2011},
which is  used in conformal field theory, renormalization and in fixing the critical spacetime dimension of string theory.
An example of zeta-function regularization is the calculation of the vacuum expectation value of the energy of a particle field in quantum field theory. More generally, the zeta-function approach can be used to regularize the full energy-momentum tensor in curved spacetime \cite{parkerQuantumFieldTheory2009}. 

On the other hand, the LD formula  \eqref{eq:LD} contains both odd- and even-order derivatives, and the latter appear with the same sign at the end-points (e.g. $t_1$ and $t_2$).  Thus,  if one attempts to construct a
composite (multi-step) formula, the even-order derivatives do not cancel out. For this reason, the LD formula is not readily\footnote{unless the function derivatives  $f'_l, f''_l,...$ can be easily expressed in terms of the function values $f_l$. } suited for converting sums to integrals. Nevertheless, the LD formula is ideally suited
for the purpose of numerical integration, i.e.
conversion of integrals to sums. In this paper,  restrict attention to  single-step, rather than composite, formulas.

For the purpose of  integration, the 
difference in accuracy between   the LD formula  \eqref{eq:LD} and the EM formula \eqref{eq:EM}  is  dramatic. Inspection of  the respective remainder terms
 \eqref{eq:LDerror} and \eqref{eq:EMerror} reveals that, 
when the function derivatives
$f^{(l)}(t)$ of all orders up to $l=n-1$ are used,
the  EM\ formula is (formally) accurate to order  $n+2$  whereas the LD formula is accurate to order
$2n$. Further inspection of the remainder term \eqref{eq:EMerror} reveals that, because the Bernoulli numbers $B_n$  decay slowly for increasing $n$, the EM formula is often only \textit{asymptotically convergent}, and in practice may fail to converge for certain functions or  large time-steps. On the other hand, the coefficients in the remainder term \eqref{eq:LDerror} decay faster than exponentially  for increasing
$n$, and  the LD formula is \textit{superconvergent}.
\begin{figure} 
        \centering
        \includegraphics[width=\linewidth]{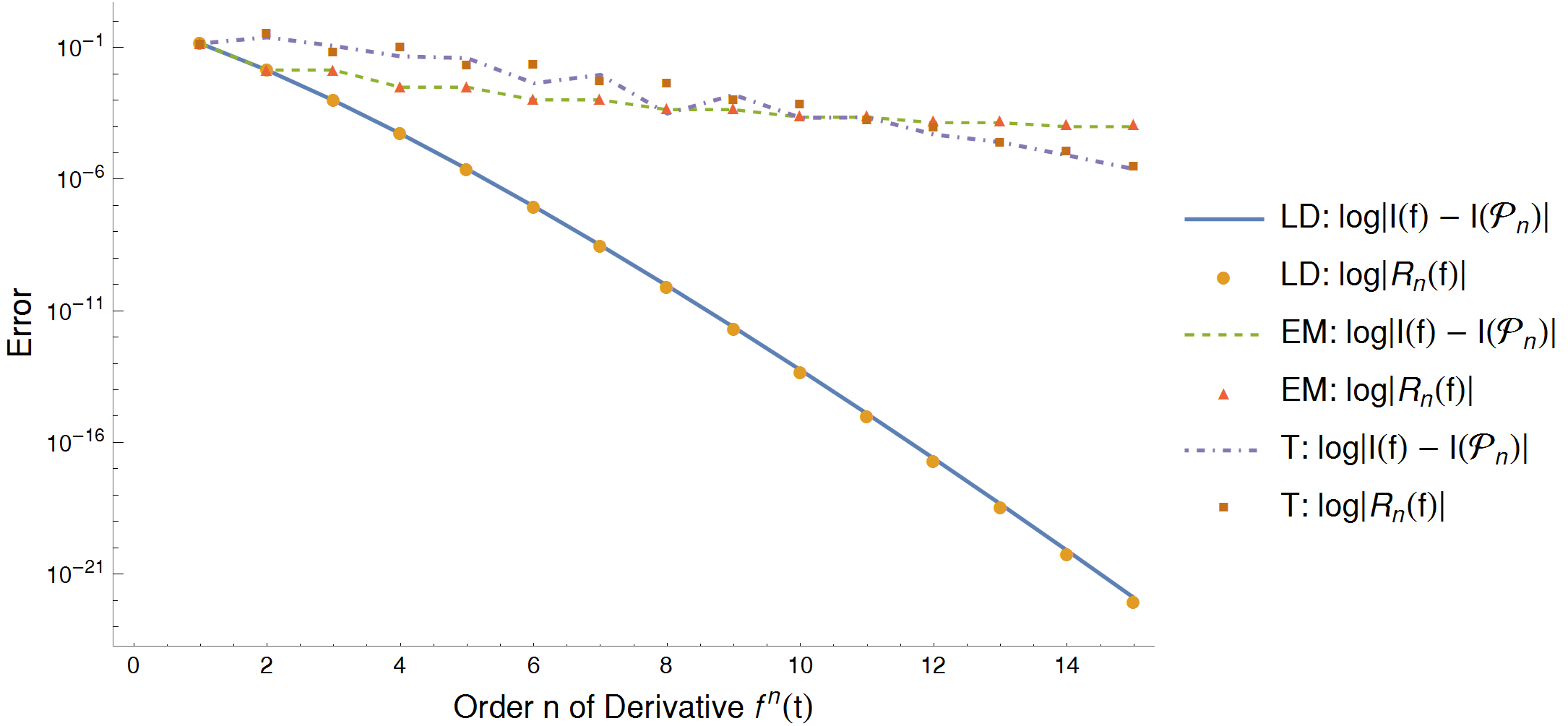}
        \caption{Illustrative comparison between  remainder terms of the Lanczos-Dyche (LD), Euler-Maclaurin (EM), Taylor Expansion (T) formulas, using derivatives of the integrand of order up to $15^{\rm th}$. 
In this example, the test function $f(t) = e^{-t^2}$ is integrated in the interval $(t_1,t_2) = (-1/2,1/2)$ with a single, large time-step $\Delta t=t_2-t_1=1.$ The solid blue, dashed green, and dashed-dot lines show the absolute difference between the exact integral  $I(f)=\int_{t_1}^{t_2} f(t) dt$ and its numerical approximation $I(\mathcal{P}_n)=\int_{t_1}^{t_2} \mathcal{P}_n (t) dt$ computed by  integrating an interpolating polynomial, that is, via the LD, EM, and T formulas, given by  Eqs.  \eqref{eq:EM},  \eqref{eq:LD} and \eqref{eq:Taylorint} respectively.  This absolute difference is in  agreement with the estimates of the respective remainder terms $R_n(f)$ for these formulas. As  the value of $\tau$  in the remainder terms makes an insignificant impact, it is chosen rather arbitrarily  to be the midpoint of the interval $\tau \in (t_{1},t_{2})$.} \label{fig:ErrComp}
\end{figure}
A comparison of the error terms of the EM and LD formulas when  computing the integral of the exponential function, shown in 
Fig. \ref{fig:ErrComp}, confirms that the latter decreases  faster than exponentially with the order $n$ of derivatives used.

\section{Ordinary differential equations}
%

\subsection{Absolute Stability}
To study stability  of a time integration scheme, we consider a linear test equation of the form \eqref{eq:dudt} with $f(t,u)=\lambda u$, $\lambda\in\mathbb{C}$, that is,
\begin{equation} \label{dudtmodel}
\frac{du}{dt} = \lambda u, \quad u(0)=1
\end{equation}
and consider the relation between $u(t_{\nu+1})$ and $u(t_{\nu})$ arising from Eqs.~\eqref{eq:ustepintf}, \eqref{eq:DI} and the relevant integration formula. For any one-step method, this relation has the form
\begin{equation} \label{eq:incru}
u(t_{\nu+1})=\zeta(\lambda \,\Delta t) u(t_{\nu}).
\end{equation}
The function $\zeta(\mu)$ is the increment   function. 
For the exact solution of Eq.~\eqref{dudtmodel}, we have $u(t_{\nu+1})=e^{\lambda \,\Delta t} u(t_{\nu})$,
and  thus
the exact increment function is\begin{equation}
\label{eq:AStabilityExact}
\zeta(\mu)=e^{\mu}.
\end{equation}
A numerical method is considered absolutely stable if the increment function for that method satisfies the condition $\left|\zeta(\lambda \,\Delta t)\right|\leq 1$.

\begin{figure}[!htpb]
        \centering
        \subfloat[\label{fig:AStabilityRK} Runge-Kutta   ]{\includegraphics[height=0.4\linewidth]{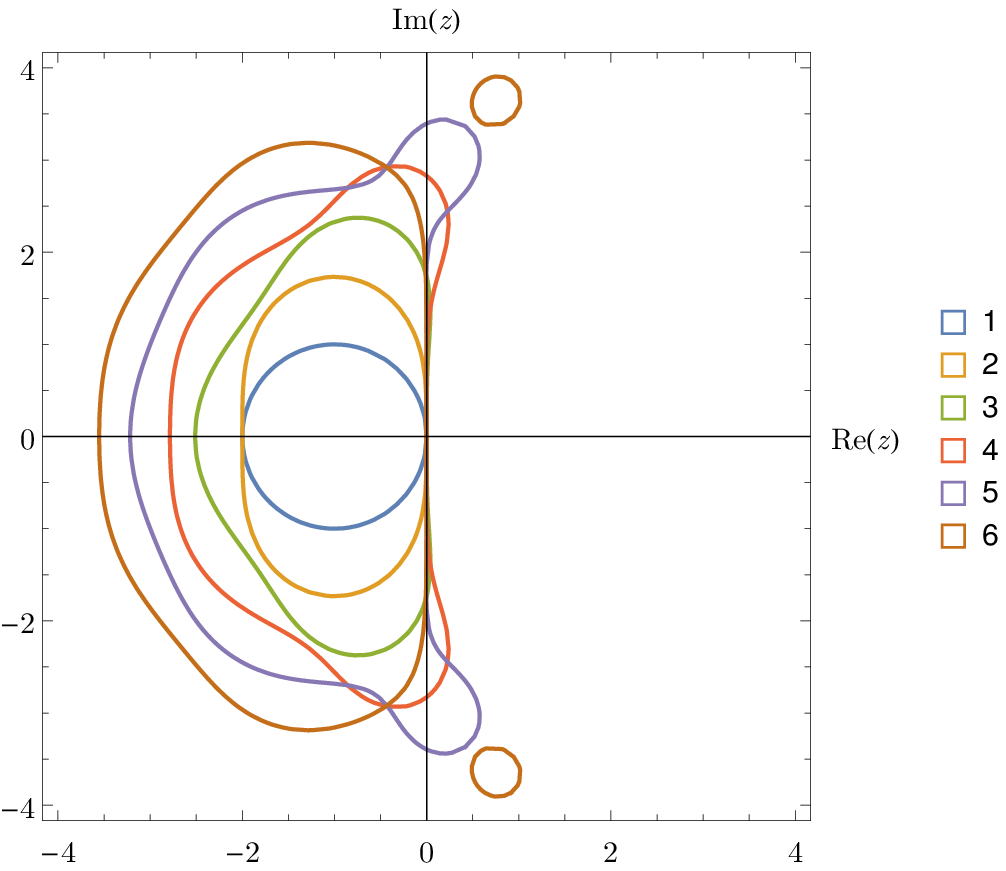}}
        \hfill
        \subfloat[\label{fig:AStabilityLD} Lanczos-Dyche
        ]{\includegraphics[height=0.4\linewidth]{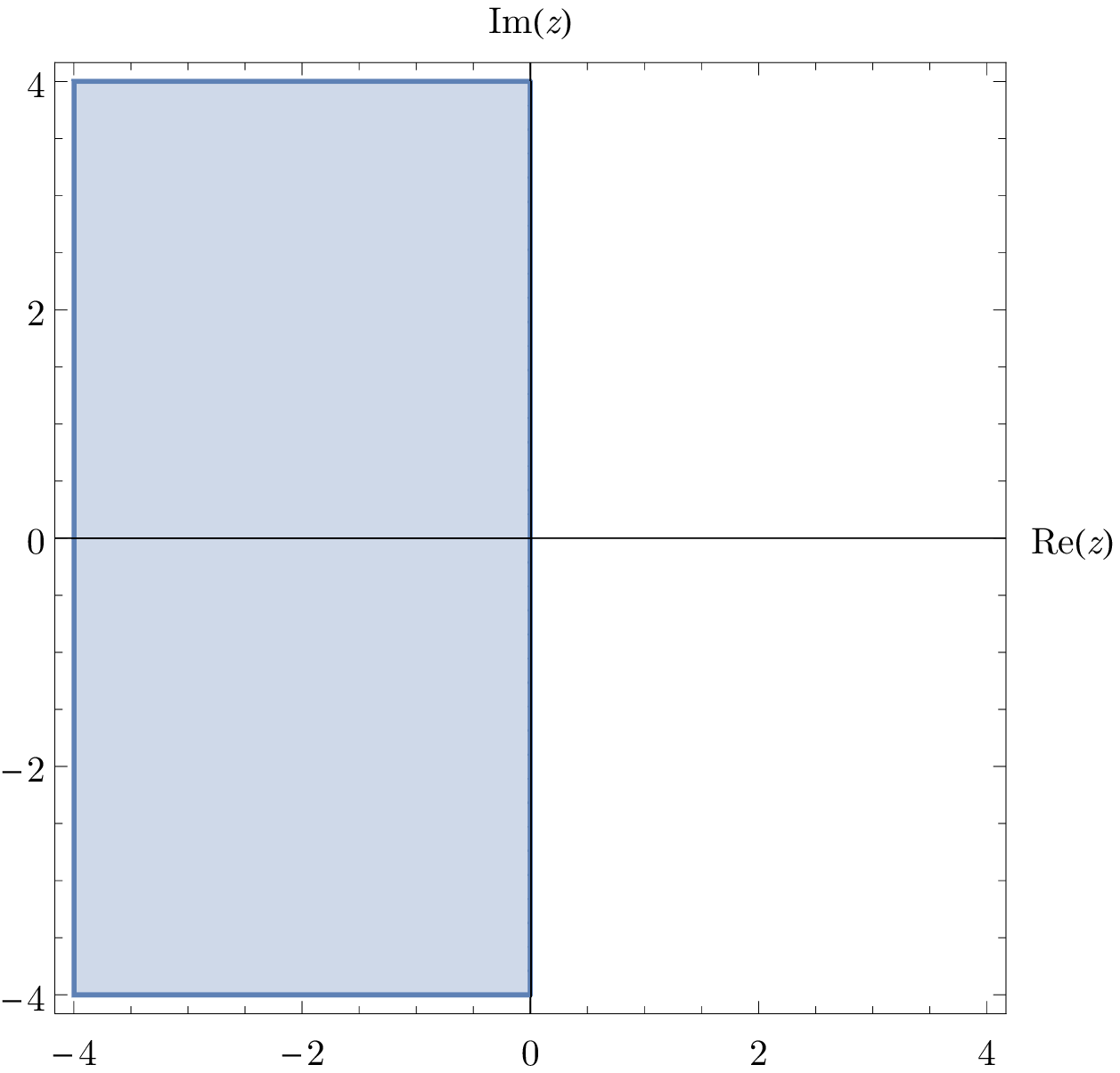}}
        \caption{\label{fig:AStability} A-stability regions for different time-stepping schemes.}
\end{figure}

For a Runge--Kutta method of order $n$, the time-stepping has the form of  Eq.~\eqref{eq:incru} with an increment function given by\begin{equation}
\label{eq:AStabilityRK}
\zeta(\mu) = \sum_{l=0}^{n}\frac{1}{l!}\mu^l,
\end{equation}
which coincides with  the Taylor expansion
of order $n$ about
$\mu =0$ of the exact increment function \eqref{eq:AStabilityExact}. 
These methods exhibit stability within the  regions plotted in Fig. \ref{fig:AStabilityRK}.

On the other hand, substituting the $l^{\rm th}$-order derivatives  $du^l/dt^{l}=\lambda^l u$ of Eq.~\eqref{dudtmodel} into the Lanczos-Dyche formula \eqref{eq:LD}, yields an increment function
\begin{equation}
\label{eq:AStabilityLD}
\zeta(\mu) = \frac{\sum_{l=0}^{n} \frac{C_{ln}}{l!}\mu^l}{\sum_{l=0}^{n} \frac{C_{ln}}{l!}\left(-\mu \right)^l},
\end{equation}
with coefficients $C_{ln}$ given by Eq.~\eqref{eq:Cln}.
The above expression coincides with   the Pad\'e approximant of order $(n,n)$ to the exact function \eqref{eq:AStabilityExact} expanded about $\mu =0$.
 It has already been proven in \cite{ehleAStableMethodsPade1973} that Pad\'e methods with stabiltity function given by \eqref{eq:AStabilityLD} are stable on the entire left half plane $\text{Re}(\mu)\leq 0$ for all orders $n\geq 1$.
Thus, the methods resulting from the Lanczos-Dyche formula \eqref{eq:LD} are A-stable (but   not L-stable, as Eq.~\eqref{eq:AStabilityLD}, does not satisfy the condition $\zeta (\mu) \rightarrow0$ as $\mu \rightarrow - \infty$). 

We note that implicit multistep methods can only be A-stable if their order is at most 2. This result, known as the second Dahlquist barrier \cite{dahlquistConvergenceStabilityNumerical1956,dahlquistStabilityErrorBounds1961,dahlquistSpecialStabilityProblem1963,butcherNumericalMethodsOrdinary2000,wannerDahlquistClassicalPapers2006},  restricts the usefulness of linear multistep methods for stiff equations. The optimal A-stable method is the trapezoidal rule, which is a special case of the Lanczos-Dyche formula for $n=1$. However, the use of higher order derivatives allows one to exceed this barrier 
\cite{brownMultiderivativeNumericalMethods1975,brownCharacteristicsImplicitMultistep1977,jeltschNoteAstabilityMultistep1976,jeltschStabilityRegionsMultistep1978,jeltschStabilityPropertiesBrown1978,jeltschA0StabilityStiffStability1979,gekelerDiscretizationMethodsStable2006}, and methods based on the Lanczos-Dyche formula are A-stable to all orders.
The property of A-stability will prove extremely useful when applying this method to partial differential equations, as it will be seen that systems discretized via the method of lines using the  Lanczos-Dyche formula are unconditionally stable, that is, these methods are not subject to a  
  Courant-Friedrichs-Lewy condition on the time-step.

We   notice that the Lanczos-Dyche formula  \eqref{eq:LD}  reduces to   Pad\'e methods  for linear differential equations. However, Pad\'e methods   rely on a rational approximation of the exponential function $e^{A t}$, and are applicable  only to linear systems. In contrast, methods resulting from the Lanczos-Dyche
formula are generally applicable to both linear and non-linear systems (albeit they are implicit in the latter case). Moreover, even for linear systems, if the system is large, Pad\'e methods provide an approximation to the  large matrix  $e^{A t}$. However, the LD formula   \eqref{eq:LD}   can be applied to each sub-equation separately, and then one can proceed to invert a linear sub-system  of much smaller size, as will be seen below in Sec.~\ref{sec:ODEs} and in more detail in a companion paper. 

\subsection{Hamiltonian systems}
Consider the Hamilton equations
\begin{subequations}
\begin{eqnarray}\label{eq:Hamiltons}
\frac{dq}{dt} \negthickspace \negthickspace&=&\negthickspace \frac{\partial H}{\partial p}
\\
\frac{dp}{dt} \negthickspace \negthickspace&=&\negthickspace \negthickspace-\frac{\partial H}{\partial q}
\end{eqnarray}
\end{subequations}
where the Hamiltonian  $H(q,p)\!:\Omega \subset \mathbb{R}^{N+1}\times  \mathbb{R}^{N+1} \rightarrow\mathbb{R}
$ is a smooth function on an open set   $\Omega$ of $  \mathbb{R}^{N+1}\times  \mathbb{R}^{N+1}
 $.

To linear order in the Taylor expansion \eqref{eq:Taylorint}, the time evolution of this system could be accomplished via  Euler integration
\begin{subequations} \label{eq:EulerHam}
\begin{eqnarray}
q_{\nu+1}\!\!&=&\!\!q_{\nu}+\Delta t \;\frac{\partial H}{\partial p_\nu}
\\
p_{\nu+1}\!\!&=&\!\!p_{\nu}-\Delta t \;\frac{\partial H}{\partial q_\nu},
\end{eqnarray}
\end{subequations}
where $\nu \in \mathbb{N}$ counts the time steps. For a typical Hamiltonian of the form
\begin{equation} \label{eq:Hamp2q2}
H=\frac{1}{2} \, p^2+V(q)
\end{equation}
the above discrete transformation is not  canonical, as the Jacobian: 
\begin{equation} \label{eq:Jacobian}
J=\frac{\partial(q_{\nu+1},p_{\nu+1})}{\partial(q_{\nu},p_{\nu})}
\end{equation}
does not preserve its unitary value, $J=1$.
Instead, the Euler integration \eqref{eq:EulerHam}  introduces
an error of order $\mathcal{O} (\Delta t^2)$ in the Jacobian \eqref{eq:Jacobian} of the transformation.


One can  include higher order terms in the Taylor expansion so the error accumulates more slowly, but symplectic structure would still not be preserved. For example, a second order Runge-Kutta integration yields the integration rule
\begin{subequations} \label{eq:RK2Ham}
\begin{eqnarray}
q_{\nu+1}\!\!&=&\!\!q_{\nu}+\Delta t \;\frac{\partial H}{\partial p_\nu} + \frac{\Delta t^2}{2} \frac{d}{dt} \bigg(\frac{\partial H}{\partial p_\nu}\bigg)
\\
p_{\nu+1}\!\!&=&\!\!p_{\nu}-\Delta t \;\frac{\partial H}{\partial q_\nu }- \frac{\Delta t^2}{2} \frac{d}{dt} \bigg(\frac{\partial H}{\partial q_\nu}\bigg),
\end{eqnarray}
\end{subequations}Since the partial derivatives of $H$ with respect to $p$ and $q$ are  functions of those variables, the chain rule:
\begin{equation}
\frac{{dH}}{{dt}} = \frac{{\partial H}}{{\partial t}} + \dot p\frac{{\partial H}}{{\partial p}} + \dot q\frac{{\partial H}}{{\partial q}}
\end{equation}
and the  equations of motion \eqref{eq:Hamiltons} may be used to remove any total time derivatives that appear at second order. A scheme of this sort would introduce an error in the Jacobian~\eqref{eq:Jacobian} of order $\Delta t^3$.

Symplectic structure may be preserved if the dependence of the Hamiltonian on the future position $q_{\nu+1}$ and momentum  $p_{\nu+1}$ at time $t_{\nu+1}$ is taken into account.  For example, one may replace the Euler method by the momentum-Verlet method \cite{apostolatosElementsTheoreticalMechanics2004}
\begin{subequations} \label{eq:VerletHamQ}
\begin{eqnarray}
q_{\nu+1}\!\!&=&\!\!q_{\nu}+\Delta t \;\frac{\partial H}{\partial p_\nu}
\\
p_{\nu+1}\!\!&=&\!\!p_{\nu}-\Delta t \;\frac{\partial H}{\partial q_{\nu+1}}.
\end{eqnarray}
\end{subequations}
or
the position-Verlet method\begin{subequations} \label{eq:VerletHamP}
\begin{eqnarray}
q_{\nu+1}\!\!&=&\!\!q_{\nu}+\Delta t \;\frac{\partial H}{\partial p_{\nu+1}}
\\
p_{\nu+1}\!\!&=&\!\!p_{\nu}-\Delta t \;\frac{\partial H}{\partial q_\nu}
\end{eqnarray}
\end{subequations}
In one of the two equations, the partial derivative of the Hamiltonian is evaluated at $t_{\nu+1}$ as opposed to $t_\nu$. For the Hamiltonian \eqref{eq:Hamp2q2}, Substituting Eq.~\eqref{eq:VerletHamQ}
or \eqref{eq:VerletHamP}  Eq.
\eqref{eq:Jacobian} yields a Jacobian $J=1$
  identically. Thus, both  Verlet methods preserve symplectic
structure. In accordance to Liouville's theorem, phase-space volume is conserved by these methods, as demonstrated in Fig.~\ref{fig:Phase}. Although these methods yield closed trajectories in phase space,  they need not conserve energy (cf. Fig.~\ref{fig:Phase}.). In addition, the resulting scheme will often be a system of nonlinear equations that must be solved numerically. We refer to the method  with advanced momentum given by Eq.~\eqref{eq:VerletHamP} as VP  (momentum-Verlet) and the method with advanced position given by Eq.~\eqref{eq:VerletHamQ} as VQ (position-Verlet) \cite{gansShadowMassRelationship2000,rowleyVariationalIntegratorsDegenerate2002,marsdenDiscreteMechanicsVariational2001,krausProjectedVariationalIntegrators2017}. 

Using the second order Lanczos-Dyche  (LD2) formula (or trapezoidal rule) to expand the solution to the Hamiltonian system \eqref{eq:Hamiltons} yields the scheme
\begin{subequations} \label{eq:LD2Ham}
\begin{eqnarray}
q_{\nu+1}\!\!&=&\!\!q_{\nu}+\frac{\Delta t}{2} \;\bigg(\frac{\partial H}{\partial p_\nu} + \frac{\partial H}{\partial p_{\nu+1}} \bigg)
\\
p_{\nu+1}\!\!&=&\!\!p_{\nu}-\frac{\Delta t}{2} \;\bigg(\frac{\partial H}{\partial q_{\nu}} + \frac{\partial H}{\partial q_{\nu+1}}\bigg).
\end{eqnarray}
\end{subequations} Note that this formula may be obtained by averaging the Verlet schemes \eqref{eq:VerletHamP} and \eqref{eq:VerletHamQ}. Intuitively, by inspecting the oval-shaped trajectories for VQ and VP in Fig.~\ref{fig:Phase}, one may expect the symmetric formula to \eqref{eq:LD2Ham} to yield circular trajectories and thus preserve symplectic structure. For linear systems, such as the  harmomic osillator, this is indeed the case. Furthermore, the formula  \eqref{eq:LD2Ham}  is time-symmetric: the  scheme does not change under the reversal  $t_\nu \leftrightarrow
t_{\nu+1}$. Since energy conservation is intimately related to time symmetry, by means of Noether's theorem, we expect this scheme to preserve energy. We theoretically and numerically explore these conjectures in the following sections.

If the fourth order LD formula (LD4) is used, one obtains
\begin{subequations} \label{eq:LD4Ham}
\begin{eqnarray}
q_{\nu+1}\!\!&=&\!\!q_{\nu}+\frac{\Delta t}{2} \;\bigg(\frac{\partial H}{\partial p_\nu} + \frac{\partial H}{\partial p_{\nu+1}} \bigg) + \frac{\Delta t^2}{12} \;\frac{d}{dt}\bigg(\frac{\partial H}{\partial p_\nu} - \frac{\partial H}{\partial p_{\nu+1}} \bigg)
\\
p_{\nu+1}\!\!&=&\!\!p_{\nu}-\frac{\Delta t}{2} \;\bigg(\frac{\partial H}{\partial q_{\nu}} + \frac{\partial H}{\partial q_{\nu+1}}\bigg) - \frac{\Delta t^2}{12} \;\frac{d}{dt}\bigg(\frac{\partial H}{\partial q_{\nu}} - \frac{\partial H}{\partial q_{\nu+1}}\bigg).
\end{eqnarray}
\end{subequations} We again expect this formula to preserve energy and symplecticity for the same reasons discussed above.

\subsection{Linear systems: harmonic oscillator} \label{sec:harmonicoscillator}






 
As a first example of a linear Hamiltonian system, we consider the simple harmonic oscillator. The Hamiltonian for this problem (with an appropriate choice of units)
is\begin{equation}
H = \frac{p^2}{2} + \frac{q^2}{2}
\end{equation} and the Hamilton equations are
\begin{subequations}
\begin{eqnarray}
\dot{q}  \!\!\! &=& \!\!\! p
\\
\dot{p} \!\!\! &=& \!\!\! -q
\end{eqnarray}
\end{subequations} Since the system is linear, the integration schemes outlined in Section (3.2) can easily be written explicitly. We explicitly   construct the first two LD schemes given by Eqs.~\eqref{eq:LD2Ham} and \eqref{eq:LD2Ham}. LD2 reads
\begin{subequations}
\begin{eqnarray}
q_{\nu+1} = q_\nu + \frac{\Delta t}{2}(p_\nu + p_{\nu+1})
\\
p_{\nu+1} = p_\nu - \frac{\Delta t}{2}(q_\nu + q_{\nu+1})
\end{eqnarray}
\end{subequations} while LD4 reads
\begin{subequations}
\begin{eqnarray}
q_{\nu+1} = q_\nu + \frac{\Delta t}{2}(p_\nu + p_{\nu+1}) - \frac{\Delta t^2}{12}(q_\nu - q_{\nu+1})
\\
p_{\nu+1} = p_\nu - \frac{\Delta t}{2}(q_\nu + q_{\nu+1}) - \frac{\Delta t^2}{12}(p_\nu - p_{\nu+1})
\end{eqnarray}
\end{subequations}Since the equations are linear, the two schemes may be explicitly solved for the advanced  momentum $p_{\nu+1}$ and  position $q_{\nu+1}$. LD2 results in
\begin{subequations}\label{eq:SHOLD2}
\begin{eqnarray}
q_{\nu+1} = q_\nu - \frac{\frac{1}{2} {\Delta t}^{2}}{1+\frac{1}{4}\Delta t^2} q_\nu+ \frac{\Delta t}{1+\frac{1}{4}\Delta t^2}p_\nu
\\
p_{\nu+1} =p_\nu - \frac{\Delta t }{1+\frac{1}{4}\Delta t^2}q_\nu- \frac{ \frac{1}{2} {\Delta t}^2~}{1+\frac{1}{4}\Delta t^2} p_\nu
\end{eqnarray}
\end{subequations} while LD4 results in
\begin{subequations}\label{eq:SHOLD4}
\begin{eqnarray}
q_{\nu+1} = q_\nu - \frac{{\frac{1}{2}\Delta t}^2}{{1+\frac{1}{12}\Delta t}^2
\left(1+\frac{1}{12}{\Delta t}^2\right) } q_\nu +\frac{{\Delta t} \left(1-\frac{1}{12}{\Delta
t}^2\right)}{{1+\frac{1}{12}\Delta t}^2 \left(1+\frac{1}{12}{\Delta t}^2\right)}p_\nu
\\
p_{\nu+1} =p_\nu  -\frac{{\Delta t} \left(1-\frac{1}{12}{\Delta
t}^2\right)}{{1+\frac{1}{12}\Delta t}^2 \left(1+\frac{1}{12}{\Delta t}^2\right)}q_\nu- \frac{{\frac{1}{2}\Delta t}^2}{{1+\frac{1}{12}\Delta t}^2
\left(1+\frac{1}{12}{\Delta t}^2\right) } p_\nu
\end{eqnarray}
\end{subequations}In both cases, polynomial expressions in $\Delta t$ were written in  Horner form, to reduce the round-off error from floating point arithmetic. Moreover, separating  out the  incremental additive change in $q$ and $p$  explicitly, also reduces round-off error from accumulated summation, and facilitates compensated summation if  desired.

The Jacobians of these transformations may be computed in the usual way, and it is found that $J=1$ for both methods. Moreover, it can be shown that 
\begin{figure} 
\centering   
\includegraphics[width=0.5\linewidth]{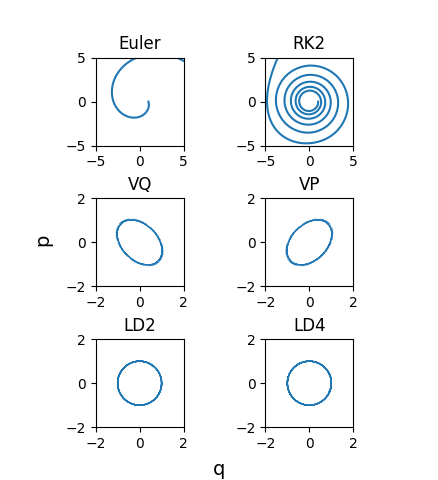}
\caption{Phase portraits of various numeric methods used to evolve the simple harmonic oscillator. Note that the Euler and RK2 methods exhibit unphysical growth. The Verlet methods create closed trajectories, but energy is not conserved throughout a period. The LD methods are symplectic and energy conserving, so they form perfect circles, the expected trajectories. To exaggerate these effects a large step size, $\Delta t = 0.75$, was used.} \label{fig:Phase}
\end{figure}

\begin{equation}
\frac{p_{\nu+1}^2}{2} + \frac{q_{\nu+1}^2}{2} = \frac{p_{\nu}^2}{2} + \frac{q_{\nu}^2}{2}
\end{equation} in both cases. Thus, the LD methods are exactly symplectic and energy conserving for this linear system. Symplecticity implies the numeric phase space trajectories should form closed curves. Energy conserving implies the Hamiltonian should have a constant value on the numeric trajectories; curves of constant energy are circles in this problem. We demonstrate this in Fig.~(\ref{fig:Phase}). The Euler (RK1) and RK2 methods are neither symplectic nor energy conserving, so their phase portraits exhibit unphysical growth. The growth in RK2 is slower than Euler because the error term is of higher order. The two Verlet methods are both symplectic, so their numeric trajectories are closed curves, but they are not energy conserving, so they are not circles as expected for this problem. The two LD methods, second and fourth order, are both symplectic and energy conserving, so their numeric trajectories are circles in phase space as expected.

\begin{figure}
 \centering
 \includegraphics[width=1.0\linewidth]{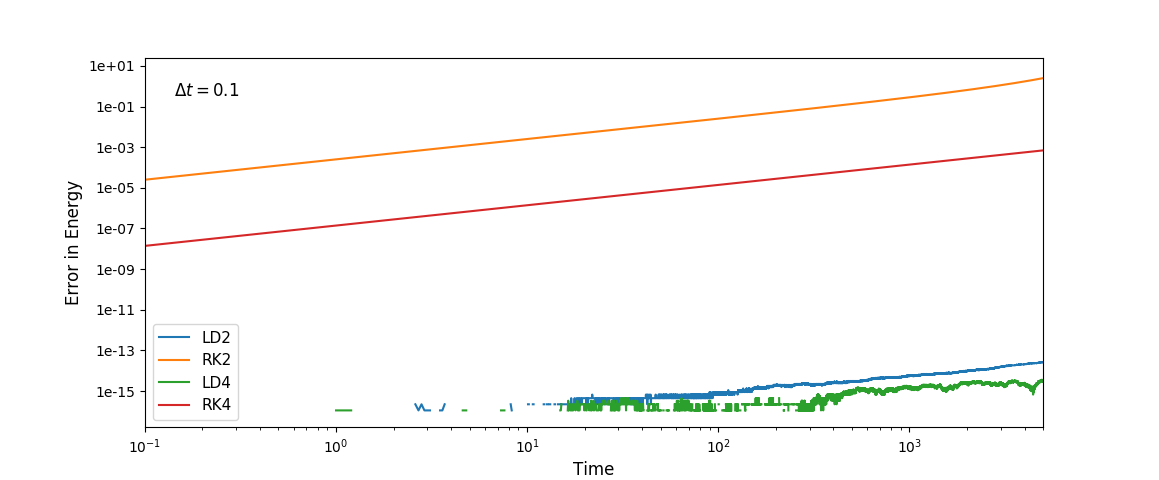}
 \caption{The relative error $\delta E /E$ in numerical energy compared with exact energy of a harmonic oscillator as a function of time. Observe that the error in the LD methods is of the order of machine precision, indicating that energy is numerically conserved. Contrast this with the RK methods, where the error grows polynomially in time, indicating that energy is not conserved.} \label{fig:SHOEnergy}
 \end{figure}

The LD methods were shown to be exactly energy conserving, so we numerically confirm this. The LD2 and LD4 methods were compared to RK2 and RK4 (the fourth order Taylor expansion). A step size of $\Delta t = 0.1$ was used to evolve the system for 5000 periods. At each time step, the relative error between the numeric energy and exact energy was computed; this quantity is defined by 
\begin{equation}
\frac{\delta E}{E} \equiv \frac{|E(t_\nu)-E(0)|}{E(0)}
\end{equation} where
\begin{equation}
E(t_\nu) \equiv \frac{p_\nu^2}{2} + \frac{q_\nu^2}{2}.
\end{equation} The RK methods were found to exhibit polynomial growth over long time scales while the LD methods conserved energy to machine precision. We note that it was necessary to use the Horner forms shown in equations \eqref{eq:SHOLD2} and \eqref{eq:SHOLD4} with compensated summation to achieve the accuracy shown in Fig.~(\ref{fig:SHOEnergy}).

 \subsection{Nonlinear systems: anharmonic oscillator}
%
%
Although the method was proved to be symplectic for linear systems, such as the harmonic oscillator, it is worth noting that this property does not generalize exactly to non-linear systems: the method is only approximately symplectic for non-linear systems. To see this, consider a general time-independent system:
\begin{equation}
\dot p =-V'(q),~~~~\dot q=p
\end{equation} 
This system is integrated to find 
\begin{equation}
p_{\nu+1} = p_{\nu} - \int_{t_{\nu}}^{t_{\nu+1}}V'(q) dt
\end{equation}
\begin{equation}
q_{\nu+1} = q_{\nu} + \int_{t_{\nu}}^{t_{\nu +1}}p~dt.
\end{equation} Let us apply the second order LD method (or trapezoidal integral rule) to approximate the above integrals by
\begin{equation}
p_{\nu+1} = p_{\nu} - \frac{\Delta t}{2} (V'(q_{\nu}) + V'(q_{\nu+1}) )
\end{equation}
\begin{equation}
q_{\nu+1} = q_{\nu} + \frac{\Delta t}{2} ( p_{\nu} + p_{\nu +1})
\end{equation}
This system can be rearranged to give equations for the advanced position $q_{\nu +1}$ and momentum $p_{\nu +1}$:
\begin{equation} \label{eq:ANCOOR}
q_{\nu+1} + \frac{\Delta t^2}{4}V'(q_{\nu+1}) = q_\nu- \frac{\Delta t^2}{4} V'(q_\nu) + \Delta t~p_\nu
\end{equation}
\begin{equation}
    p_{\nu+1} = p_{\nu} - \frac{\Delta t}{2}(V'(q_{\nu}) + V'(q_{\nu+1}))
\end{equation}

The Jacobian for this method can be calculated as
\begin{equation}\label{eq:JAC}
J = \frac{1+\frac{\Delta t^2}{4}V''(q_{\nu})}{1+\frac{\Delta t^2}{4}V''(q_{\nu+1})}
\end{equation}We see that the method is not symplectic unless $V''(q_{\nu})=V''(q_{\nu+1})$,
which is only true if $V$ is either a constant, linear, or quadratic function of position. 

Nevertheless, we emphasize that the method is still superior to explicit methods. Consider the pendulum equation, with  potential $V(q)=\cos q$, as an example of an anharmonic oscillator. If we assume $\Delta t$ is small, then Eq.~\eqref{eq:JAC} may be expanded as
\begin{equation}\label{eq:JACEX}
 J \simeq 1 - \frac{\Delta t^2}{4}(\cos q_{\nu+1} - \cos q_{\nu}) = 1+ \mathcal{O}(\Delta t^3)
\end{equation}
Although the order of  deviation from symplecticity is formally the same as the order of the method,  the method is  more accurate for this problem than explicit methods, such as the ordinary Runge-Kutta methods. This is because of two reasons:\ (i) For the LD methods, the error in the Jacobian is oscillatory
and thus has an upper bound, in contrast to the RK methods, where the error grows in time without bound. (ii) As seen by inspecting the  remainder terms of each formula, the error for the LD methods is  significantly lower compared to RK methods of the same order. The difference increases by several orders of magnitude at high order, as demonstrated in the previous section.\   

We demonstrate this explicitly by comparing the second and fourth order LD methods to second and fourth order RK methods. To carry out the LD integration, equation \eqref{eq:ANCOOR} and the corresponding fourth order equation needed to be solved numerically. A simple Newton-Raphson algorithm was adequate. For the second order method we iterated the equality
\begin{equation}
q_{i+1} = q_i - \frac{q_i + \frac{\Delta t^2}{4} \sin q_i - (q_\nu + \Delta t ~p_{\nu} +\frac{\Delta t^2}{4} \sin q_{\nu})}{1+\frac{\Delta t^2}{4} \cos q_i}
\end{equation}
with $q_i$ initialized to $q_{\nu}$. We found that the algorithm converged rapidly, within three iterations, regardless of $\Delta t$'s magnitude. A similar procedure was followed for the fourth order LD method. 

\begin{figure}
\centering
\includegraphics[width=\linewidth]{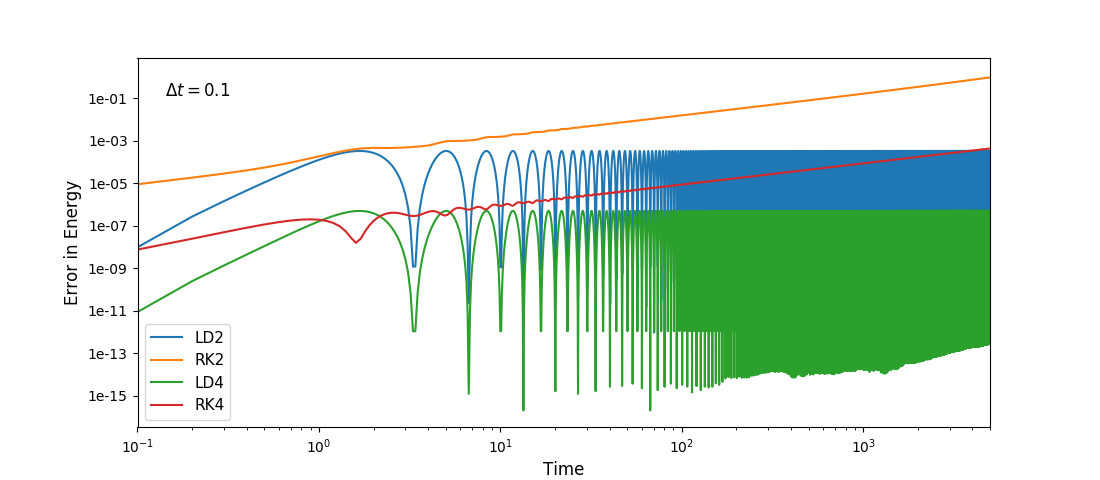}
\caption{The relative error  $\delta
E/E$  in the numerical energy compared to the exact energy as a function of time for the pendulum equation. In all simulations $\Delta t = 0.1$. Note that for the LD methods the signal is locally a periodic function, strongly suggesting that the method is effectively symplectic, even if not exactly symplectic.  Note that that the magnitude of the relative error is small. Contrast this with the RK methods which quickly display polynomial growth.}
\label{fig:PEN}
\end{figure}

We plot the deviation from the true energy as a function of time from a numeric simulation in Fig.~(\ref{fig:PEN}). We note that the deviation from the true energy is a periodic function of time for both LD methods and remains bounded, as opposed to  explicit methods for which the deviation grows polynomially in time. This demonstrates that our LD methods are still superior to explicit methods for the pendulum equation, despite their non-exact symplecticity.

This can be understood by inspecting Eq.~\eqref{eq:JACEX}, which suggests  that the deviation from symplecticity is formally of order $\mathcal{O}(\Delta t^3)$, but oscillates sinusoidally with $q_{\nu}$ and  $q_{\nu+1}$. Each of these contributions  averages to zero over one period, so the second order LD method can be considered ``symplectic on average''. A similar result holds for the fourth order LD method, since the  generalized force is sinusoidal.

\subsection{Nonconservative Systems: Damped Oscillator}
Having demonstrated the method's utility for both linear and nonlinear conservative systems, we now turn to nonconservative systems. Although there is no longer a manifestly conserved quantity, such systems can still be Hamiltonian and thus preserve symplectic structure. We expect our method to more useful in these problems than an explicit method, such as Runge-Kutta, for this reason.

Consider the damped harmonic oscillator with unit mass and spring constant as a specific example:
\begin{equation}
\ddot{q} + \gamma \dot{q} + q = 0.
\end{equation}
The damping force cannot be expressed as the derivative of a potential, so this system does not conserve mechanical energy. However, one can still find a Lagrangian that gives rise to this equation
\cite{havasRangeApplicationLagrange1957,leachDirectMethodDetermination1977,lemosCanonicalApproachDamped1979,leubnerInequivalentLagrangiansConstants1981,apostolatosElementsTheoreticalMechanics2004,mcdonaldDampedOscillatorHamiltonian2015,choudhuriSymmetriesConservationLaws2008}:
\begin{equation}
L = e^{-\gamma t}\bigg(\frac{\dot{q}^2}{2} - \frac{q^2}{2}\bigg).
\end{equation}
One can then compute a Hamiltonian,
\begin{equation}
H = \frac{p^2}{2}e^{-\gamma t} + \frac{q^2}{2}e^{\gamma t},
\end{equation}
and canonical equations of motion,
\begin{equation}
\dot{p} = -e^{\gamma t}q, ~~~\dot{q} = e^{-\gamma t}p.
\end{equation}
Note that $p$ is  the canonical momentum, distinct from the kinematic momentum $\dot q$ in this problem. Although energy is not conserved, one can show that the following quantity is a constant of motion \cite{qureshiExactEquationMotion2010}:
\begin{equation}
C = \frac{p^2}{2}e^{-\gamma t} + \frac{\gamma}{2}pq + \frac{q^2}{2}e^{\gamma t}.
\end{equation}

\begin{figure}
\centering
\includegraphics[width=\linewidth]{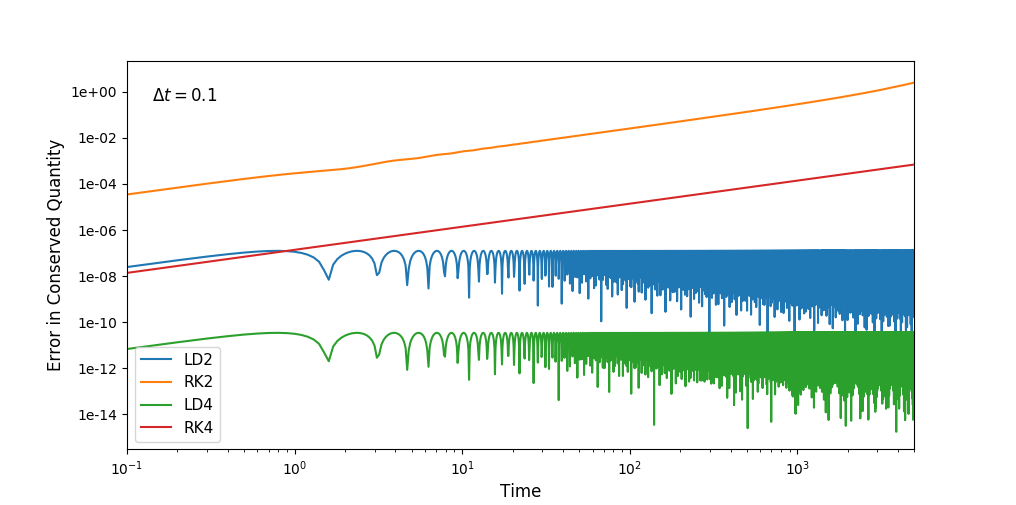}
\caption{The relative error $\delta C/C$ in the numerical value of the conserved quantity $C$, compared to its exact value, as a function of simulation time for the damped oscillator. For all simulations, $\gamma = 10^{-4}$ and $\Delta t = 0.1$. As with the pendulum equation, the error in the LD methods is bounded and oscillatory, indicating effective symplecticity (even if it is not exactly symplectic) while the RK methods exhibit polynomial error growth.} \label{fig:DAMP}
\end{figure}

Since this system is Hamiltonian, its symplectic structure should be preserved, and the constant $C$ should be conserved through the  evolution. This provides a criterion to test numeric schemes. In Figure (\ref{fig:DAMP}) we plot the relative error in $C$ as a function of time resulting from second and fourth order LD and RK methods. Like the pendulum equation, the LD scheme for this problem is not exactly conservative or symplectic, so $C$ exhibits oscillations. As mentioned earlier, the oscillations are bounded, and increasing the method's order reduces the amplitude of these oscillations. Contrast this with the RK methods, in which $C$ grows polynomially in time.

The performance of the LD methods may be compared to that of other variational integrators, such as the `slimplectic' integrator \cite{tsangSlimplecticIntegratorsVariational2015}.
For dissipative linear systems, such as the damped harmonic oscillator,  numerical accuracy is essentially the same. For linear conservative systems, the LD methods conserve both energy and symplectic structure (i.e. the Jacobian \ref{eq:Jacobian} of the canonical transformation associated with motion)
exactly, while  most symplectic integrators typically conserve the Jacobian exactly but the energy is conserved approximately (in the sense that its growth is bound) \cite{marsdenDiscreteMechanicsVariational2001}. On the other hand, for non-linear systems,  LD methods conserve energy and the Jacobian  approximately,  but symplectic methods typically conserve the Jacobian exactly.   A discussion of the connection between symmetric and symplectic integration may be found in \cite{hairerSymplecticIntegrationHamiltonian2002}.

\section{Systems of  differential equations}

\subsection{Ordinary differential equations: coupled harmonic oscillators} \label{sec:ODEs}
The properties demonstrated above generalize to systems of ordinary of partial differential equations.
For example, let us consider a 1-dimensional system of coupled harmonic oscillators. The Hamiltonian can be written in matrix  form
as\begin{equation}
H(q,p)= \frac{1}{2} (M^{-1})^{ab} p_a p_b + \frac{1}{2} K_{ab}q^a q^b
\end{equation}
where $\mathbf{M}$  and $\mathbf{K}$ are constant matrices that depend on the mass and Hooke constant of each oscillator, and the indices $a,b,...$ are summed over the different oscillators. The Hamilton equations of motion,
\begin{subequations} \label{eq:coupledharmonicoscillator}
\begin{eqnarray}
\dot{q}_a  \!\!\! &=& \!\!\! (M^{-1})^{ab} p_b
\\
\dot{p}_a \!\!\! &=& \!\!\! -K_{ab}q^b,
\end{eqnarray}
\end{subequations}
can be integrated in time to any order using the LD formula \eqref{eq:LD},
and the equations of motion can be used to eliminate all time derivatives.
This procedure, after solving the resulting linear system, yields an explicit scheme
\begin{subequations}
\begin{eqnarray}
\mathbf{q}_{\nu+1}  \!\!\! &=& \!\!\! \mathbf{A} \, \mathbf{q}_{\nu} + \mathbf{B} \, \mathbf{p}_{\nu}
\\
\mathbf{p}_{\nu+1}  \!\!\! &=& \!\!\! \mathbf{C} \, \mathbf{q}_{\nu} + \mathbf{D} \, \mathbf{p}_{\nu},
\end{eqnarray}
\end{subequations}
where the matrices $\mathbf{A},...,\mathbf{D}$ can be computed using the methods outlined in Sec. \ref{sec:harmonicoscillator}. 
This example serves as a prelude to the 
symplectic integration of hyperbolic PDEs such as the wave equation. For example, upon discretization via the method of lines (cf. Sec. \ref{sec:MOL}), the Klein-Gordon equation takes  the form of a  system of coupled ODEs identical to Eq.~\eqref{eq:coupledharmonicoscillator}.

\subsection{Partial differential equations: method of lines} \label{sec:MOL}
As mentioned above, a common approach to solving systems of evolution equations is the
method of lines. For concreteness, let us consider a first-order in time
scalar differential equation in 1+1 dimensions
\begin{equation} \label{eq:GENERAL}
\frac{\partial}{\partial t}
 u(t,x)= L (u(t,x))
\end{equation}
where $ L$ is a (not necessarily linear) spatial differential operator. If
one uses finite-difference or pseudospectral methods for the spatial derivatives
\cite{bengtfornbergPracticalGuidePseudospectral1996,markakisHighorderDifferencePseudospectral2014}, one may rewrite the original partial differential
equation
as a coupled system of ordinary differential equations of the form

\begin{equation}
\frac{d \mathbf{u}(t)}{d t}
=  \bold L (\bold u(t))
\end{equation}

\begin{equation} \label{eq:dudtisLu}
\frac{d u_i(t)}{d t}
=  \sum_{j=0}^{N} L_{ij} (u_j(t))
\end{equation}
where $u_i(t)=u(t,x_i)$
is a vector constructed from the values of the scalar field $u(t,x)$ in the
spatial grid points $x_i$ and $\bold L$ is a matrix (band-diagonal in the
case of finite differences, or full in the case of pseudo-spectral methods)
that couples the different grid points.

Eq.~\eqref{eq:dudtisLu} can be time-stepped using the fundamental theorem of calculus, in accordance with Eq.~\eqref{eq:ustepintf}.
Thus, the methods applied in the previous sections to solve ordinary differential equations, can be generalized to partial differential equations. A detailed derivation and demonstration of the symmetry, symplecticity and stability properties afforded by the LD methods in such PDEs will be given in a companion paper. It will be shown that the method is particularly powerful for linear or quasi-linear PDEs, as it is unconditionally stable, and thus overcomes the Courant-Friedrichs-Lewy
limit, and moreover preserves properties of the continuum system, such as energy or probability in  the Schr\"odinger equation and other parabolic or hyperbolic partial differential equations. 

\section{Conclusions}

To our knowledge, although 
Pad\'e approximation methods have been known for linear PDEs, their derivation
from the LD formula has not been brought forward before. More
importantly, since 
Pad\'e methods amount to a 
Pad\'e expansion of the exponential function, their applicability is restricted
to  linear systems. In contrast, the LD  formula only assumes that the problem has continuous derivatives up to some order, so it
is equally valid for non-linear systems.
Moreover, even in the context of linear systems, the LD methods are more general than Pad\'e methods, as the former can be straightforwardly applied  to each sub-equation separately, requiring inversion of matrices of smaller size. This is not straightforward with Pad\'e methods.
This property is particularly useful for integrating periodic Hamiltonian systems or wave-like equations (cf. Sec.~\ref{sec:ODEs}), as will be demonstrated in a companion paper.

In the context of ODEs, the LD methods are examples of  ``multi-derivative methods'', as the schemes depend on derivatives of the function in question. The remarkable stability properties of such methods have been studied by Brown et al. \cite{brownMultiderivativeNumericalMethods1975,brownCharacteristicsImplicitMultistep1977,jeltschStabilityPropertiesBrown1978,jeltschA0StabilityStiffStability1979}. But in addition, their time-symmetric nature gives them great promise for symplectic and energy-conserving integration. 
Our hope is to draw attention to the properties of these schemes by demonstrating their utility for integrating Hamiltonian systems and  many other problems of physical interest. We have shown that the LD methods are exactly  symplectic and energy conserving for linear problems and outperform explicit schemes for nonlinear and dissipative problems. Their remarkable accuracy  for quadratic Hamiltonians makes them  very useful for systems of linear or quasi-linear ODEs.

Moreover, when  linear PDEs are discretized  using the method of lines, a system of coupled linear ODEs is obtained. Many such problems are endowed with symplectic structure analogous to that outlined in Sec.~\ref{sec:ODEs}. This suggests that an integration scheme with the same advantages as for ODEs (stability, symplecticity, and energy conservation) is often possible. In a companion paper, we will use the LD  formula to obtain such schemes for  prototypical hyperbolic or parabolic PDEs, such as the advection, wave, and  Schr\"odinger equations.
While LD methods certainly posses utility for ODEs, we will demonstrate that it is in the numerical integration of PDEs that these schemes are most powerful.

Finally, we emphasize that, although much of the discussion has focused on linear problems, the LD formula is equally valid nonlinear systems. It will be interesting to explore the use of such integration schemes in nonlinear partial differential equations.






\section*{Acknowledgements}
We thank Theocharis Apostolatos, Petros J. Ioannou, Leo Stein and George Pappas for valuable discussions.
C.M. was supported by the European Union's Horizon 2020 research and innovation programme under the Marie Sk\l odowska-Curie grant agreement No 753115.
H.S. was supported in part by the Jade Mountain Young Scholar Award No 107V0201.

\section*{References}

\bibliographystyle{elsarticle-num}
\bibliography{numerical}

\end{document}